\documentclass[letterpaper,11pt]{article}
\usepackage{theorem}
\usepackage{epsfig}
\usepackage{amsmath}
\usepackage{amssymb}
\usepackage{amsfonts}
\usepackage{pst-node}

\newtheorem{thm}{\bfseries Theorem}
\newtheorem{lem}[thm]{\bfseries Lemma}        
\newtheorem{cor}[thm]{\bfseries Corollary}     

\newtheorem{cl}[thm]{\bfseries Claim}

\def\cl{\mathop{\mathrm{cl}}\nolimits}

\newenvironment{proof}{\begin{rm}\par\medskip          
\noindent{\scshape Proof:}}{\quad $\Box$\end{rm}\medskip}  
\newenvironment{ack}
{\begin{rm}\par\bigskip\noindent{\bf Acknowledgements.}\quad\em}{\end{rm}}

\title{
A linear equation for Minkowski sums of polytopes relatively in general
position%
\footnote{This research was supported by the Swiss National Science Foundation Project 200021-105202,
``Polytopes, Matroids and Polynomial Systems''.}
}
\date{November 25, 2007, Revised April 28, 2008}
\author{Komei Fukuda%
\footnote{Also affiliated with Institute for Operations Research
and Institute of Theoretical Computer Science
ETH Zentrum, Zurich, Switzerland.}
\\
Mathematics Institute\\
EPFL, Lausanne\\Switzerland\\
komei.fukuda@epfl.ch\\
\and 
Christophe Weibel\\
Department of Mathematics\\
McGill University, Montreal\\Canada\\
weibel@math.mcgill.ca
}

\begin{document}
\maketitle

\begin{abstract}
The objective of this paper is to study a special family of Minkowski
sums, that is of polytopes relatively in general position. We show
that the maximum number of faces in the sum can be attained by this
family. We present a new linear equation that is satisfied by
$f$-vectors of the sum and the summands. We study some of the
implications of this equation.
\end{abstract}
\section{Introduction}
Minkowski sums of polytopes naturally arise in many domains, ranging
from mechanical engineering \cite{Petit04} and robotics
\cite{lozano79} to algebra \cite{Gritzmann93,Sturmfels98}.  In direct
applications to physical models, the important factor is often the
general shape of the sum, and some details can be approximated. In
theory applications, though, it is often the combinatorial structure
of the sum which is relevant.

The scope of this paper centers on that combinatorial structure.  Few
results are known as yet as to the structure of Minkowski sums.  It is
usually difficult to estimate even the number of $k$-dimensional faces
($k$-faces) of the sum, let alone the general structure. This paper
focuses on a certain family of sums and its properties, which can be
used to make general statements about Minkowski sums.

Every nonempty face of a Minkowski sum can be \emph{decomposed} uniquely
into a sum of faces of each summand \cite{Fukuda04}. We say this
decomposition is \emph{exact} when the dimension of the sum is equal
to the sum of the dimensions of the summands. When all facets have an
exact decomposition, we say the summands are \emph{relatively in
general position}.

Our first observation is that the maximal number of faces in a sum can
be attained when summands are relatively in general position.

For any polytope $P$, we denote its \emph{f-vector} by $f(P)$ whose
$k$th component $f_k(P)$ is
the number of faces of dimension $k$ in $P$ for $k=-1,0,\ldots,\dim(P)$.
For any set $S\subseteq
\{-1,\ldots,\dim(P)\}$, we define
$f_S(P)$ as the number of chains in $P$ in which the dimensions of
elements of the chain are exactly the elements of $S$.
The vector with components $f_S(P)$ for all $S\subseteq
\{-1,\ldots,\dim(P)\}$ is called the \emph{extended f-vector}
and denoted also as $f(P)$. 

\begin{thm}\label{maxthm}
Let $P=P_1+\cdots+P_r$ be a Minkowski sum. There is a Minkowski sum
$P'=P_1'+\cdots+P_r'$ of polytopes relatively in general position so
that $f_k(P_i')=f_k(P_i)$ for all $i$ and $k$, and so that
$f_k(P')\geq f_k(P)$ for all $k$.
\end{thm}
This family can therefore be used for computing the maximum complexity
of Minkowski sums. We present now our main theorem.

For any sum of Minkowski $P=P_1+\cdots+P_r$, and for any face $F$ of
$P$, we will define its $f^\delta$-vector as
$$
f^\delta_k(F)=f_k(F)-(f_k(F_1)+\cdots+f_k(F_r)),
$$
where $F_1,\ldots,F_r$ is the decomposition of $F$ in faces of
$P_1,\ldots,P_r$ respectively.

\begin{thm}\label{mainthm}
Let $P_1,\ldots,P_r$ be $d$-dimensional polytopes relatively in general
position, and $P=P_1+\cdots+P_r$ their Minkowski sum. Then
$$
\sum_{k=0}^{d-1}(-1)^k kf^\delta_k(P) = 0.
$$
\end{thm}
Note that the form is rather similar to Euler's Equation:
$$
\sum_{k=0}^{d-1}(-1)^k f_k(P) = 1 - (-1)^d.
$$
By using Euler's Equation, we can write the theorem slightly differently.
\begin{cor}\label{maincor}
Let $P_1,\ldots,P_r$ be $d$-dimensional polytopes relatively in general position, and
$P=P_1+\cdots+P_r$ their Minkowski sum. Then for all $a$,
$$
\sum_{k=0}^{d-1}(-1)^k (k+a)f^\delta_k(P) = a(1-r)(1 - (-1)^d).
$$
\end{cor}

Additionally, the main theorem can be extended to sums of polytopes
which are not full-dimensional:

\begin{thm}\label{maincor2}
Let $P_1,\ldots,P_r$ be polytopes relatively in general position, and $P=P_1+\cdots+P_r$
their $d$-dimensional Minkowski sum. Furthermore, let $S\subseteq\{1,\dots,r\}$
be the set of indices $i$ for which $\dim(P_i)<d$. Then
$$
\sum_{k=0}^{d-1}(-1)^k k f^\delta_k(P) = (-1)^{d+1}\sum_{i\in S}\dim(P_i).
$$
\end{thm}

The main theorem has an interesting application when used in
conjunction with the following theorem about perfectly centered
polytopes (defined in the next section). For any face $F$ of a polytope
$P$ containing the origin, we denote as $F^D$ the associated dual face
of the dual polytope $P^*$.

\begin{thm}[\cite{Fukuda07}]\label{thm:oldthm}
Let $P$ be a perfectly centered polytope.  A subset $H$ of $P+P^*$ is
a nontrivial face of $P+P^*$ if and only if $H=G+F^D$ for some ordered
nontrivial faces $G\subseteq F$ of $P$.
\end{thm}
A perfectly centered polytope and its dual satisfy the general
position condition posed by the main theorem, which makes it a
statement about lattices of polytopes combinatorially similar to a
perfectly centered polytope. We extend the statement to all polytopes
and Eulerian posets:

\begin{thm}\label{nestthm}
Let $f$ be the extended $f$-vector of an Eulerian poset of rank $d$. Then
$$
\sum_{k=0}^{d-1} (-1)^k k \left(\sum_{i=0}^kf_{i,i+d-1-k}-f_k-f_{d-1-k}\right) = 0.
$$
\end{thm}
\section{Normal cones}\label{fans}
For each nonempty face $F$ of a convex polytope $P$, we define its
\emph{outer normal cone} $\mathcal{N}(F;P)$ as the set of vectors
defining linear functions which are maximized over $P$ on $F$. If $P$
has full dimension, the dimension of $\mathcal{N}(F;P)$ is
$\dim(P)-\dim(F)$.

We say a polytope $P$ is \emph{perfectly centered} if for any nonempty
face $F$ of $P$, the intersection $relint(F)\cap\mathcal{N}(F;P)$ is
nonempty.  We will use this definition in
Section~\ref{sec:ext}. Perfectly centered polytopes are further
studied in \cite{Fukuda07}.

The set of all normal cones of a polytope $P$ is a polyhedral complex
of relatively open cones whose body is $\mathbb{R}^d$, and which is
known as the \emph{normal fan} of $P$. The combinatorial structure
(i.e. the face poset) of the normal fan is dual to that of the
polytope, excluding the empty face.

It is not difficult to see that the normal fan of a Minkowski sum is
the common refinement of normal fans of the summands, i.e. the set of
nonempty intersections of normal cones of the faces of the summands
(See e.g. Section~7.2 of \cite{Ziegler95}).  For any face $F$ of a
Minkowski sum $P=P_1+\cdots+P_r$, $F$ decomposes into $F_1+\cdots+F_r$
if and only if
$$
\mathcal{N}(F;P) = \mathcal{N}(F_1;P_1)\cap \cdots\cap \mathcal{N}(F_r;P_r).
$$

\section{Maximality}\label{max}
We now show that the maximum number of faces of a Minkowski sum can
always be attained when the summands are relatively in general
position.

\begin{quote}
\textbf{Theorem~\ref{maxthm}}~
\emph{Let $P=P_1+\cdots+P_r$ be a Minkowski
sum. There is a Minkowski sum $P'=P_1'+\cdots+P_r'$ of polytopes
relatively in general position so that $f_k(P_i')=f_k(P_i)$ for all
$i$ and $k$, and so that $f_k(P')\geq f_k(P)$ for all $k$.
}
\end{quote}
\begin{proof}
Let $P=P_1+P_2$ be a Minkowski sum of polytopes not relatively in
general position. We show that if we rotate $F_1$ by a small angle on
an axis in general position, then the number of faces won't diminish.

Let $F=F_1+F_2$ be a face whose decomposition is not exact, that is
$\dim(F)<\dim(F_1)+\dim(F_2)$. In terms of normal cones, it means that
$\mathcal{N}(F;P)= \mathcal{N}(F_1;P_1)\cap \mathcal{N}(F_2;P_2)$, and
$\dim(\mathcal{N}(F;P))>
\dim(\mathcal{N}(F_1;P_1))+\dim(\mathcal{N}(F_2;P_2))-d$.  If we perturb
$P_1$ by a small enough rotation on an axe in general position, there
is a superface $\mathcal{N}(G_2;P_2)$ of $\mathcal{N}(F_2;P_2)$ so
that $\mathcal{N}(F_1;P_1)\cap \mathcal{N}(G_2;P_2) \ne \emptyset$,
and so that $\dim(\mathcal{N}(F_1;P_1)\cap \mathcal{N}(G_2;P_2))=
\dim(\mathcal{N}(F;P))$. This means $F_1$ and $G_2$ sum to a face $F'$
with $\dim(F')=\dim(F)$. So for every face with an inexact
decomposition, there is a new face of the same dimension with an exact
decomposition.  If the angle is small enough, every face with an exact
decomposition still exists. Therefore the number of faces won't
diminish, and the new sum is relatively in general position.

By induction, we can slightly rotate the summands $P_1,\ldots,P_r$ so
that they are relatively in general position, without diminishing the
number of faces in their sum.
\end{proof}
\section{Proof}\label{proof}
We prove in this section the main theorem, first when all summands are
full-dimensional, then extending it to the general case.

We start by an important lemma about the $f^\delta$-vector of faces of
a Minkowski sum which have an exact decomposition:
\begin{lem}\label{lemface}
Let $F=F_1+\cdots+F_r$ be a nonempty face of the Minkowski sum
$P=P_1+\cdots+P_r$ with an exact decomposition. Then:
$$
\sum_{k=0}^{d-1}(-1)^k kf^\delta_k(F) = 0.
$$
\begin{proof}
Since $F$ has an exact decomposition, all of its subfaces also have
one. Furthermore, for any set $(G_1,\ldots,G_r)$ so that $G_i\subseteq
F_i$ for all $i$, the sum $G=G_1+\cdots+G_r$ is a subface of $F$.

Let $d_i$ be the dimension of $F_i$ and $f^i$ its $f$-vector.  It can
be written as $(f^i_0,\ldots,f^i_{d_i})$, with $f^i_{d_i}=1$.  The
$f$-vector $f^i$ verifies Euler's equation, which means
$$
\sum_{k=0}^{d_i}(-1)^{k}f^i_k=1.
$$
Let us define the characteristic function $p_i(x)$ of the vectors
$f^i$ as follows:
$$
p_i(x) = f^i_0x^0+\cdots+f^i_{d_i}x^{d_i}.
$$
Euler's equation can now be written as $p_i(-1)=1$.

Let $f$ be the $f$-vector of $F$. Since any $r$-tuple of subfaces
$G_i$ of $F_i$ sums to a subface of $F$, we can write:
$$
f_i = \sum_{e_1+\cdots+e_r=i}(f^1_{e_1}\cdots f^r_{e_r}).
$$
Therefore, if we denote as $p(x)$ the characteristic function of the
$f$-vector of $F$, we have that $p(x)=\prod_{i=1}^r p_i(x)$.

If we denote as $p_\delta(x)$ the characteristic function of the
$f^\delta$-vector of $F$, we have that
$p_\delta(x)=f^\delta_0x^0+\cdots+f^\delta_{\dim(F)}x^{\dim(F)}=p(x)-(p_1(x)+\cdots+p_r(x))=
\prod_{i=1}^r p_i(x)-\sum_{i=1}^r p_i(x)$.  It is easy to see that
$$
\sum_{k=0}^{\dim(F)}(-1)^k kf^\delta_k(F) = -(p_\delta)'(-1).
$$
Since $p'(x)=\sum_{i=1}^r \left(p_i'(x) \prod_{j\neq i}p_j(x)\right)$,
and we have $p_j(-1)=1$ for all $j$,
$p_\delta'(-1)=\sum_{i=1}^r p_i'(-1)-\sum_{i=1}^r p_i'(-1)=0$.
\end{proof}
\end{lem}
As we can see, the $f^\delta$-vector of the faces of a Minkowski sum
which have an exact decomposition are so to say transparent to the
equation of the final theorem. We will now show that the
$f^\delta$-vector of a Minkowski sum of $d$-polytopes can be written
as an alternated sum of the $f^\delta$-vector of its proper faces.

\begin{thm}\label{thm:delta}
Let $P_1,\ldots,P_r$ be $d$-dimensional polytopes and
$P=P_1+\cdots+P_r$ their Minkowski sum. Then for any $k<d$,
$$
f^\delta_k(P)=-\sum_{F\subset P}(-1)^{d-\dim(F)}f^\delta_k(F).
$$
Or equivalently, we can write
$$
\sum_{F\subseteq P}(-1)^{d-\dim(F)}f^\delta_k(F)=0.
$$
\end{thm}

The $f^\delta$-vector of $P$ consists of the difference between its
$f$-vector and that of its summands.  The equation above actually holds for
each of these $f$-vectors. We are going to prove them separately in
the two following lemmas.

First we prove that Theorem~\ref{thm:delta} holds 
for the $f$-vector of the sum:
\begin{lem}\label{lem:euler}
Let $P$ be a $d$-dimensional polytope. Then, for any $k<d$,
$$
\sum_{F\subseteq P}(-1)^{d-\dim(F)}f_k(F)=0.
$$
\begin{proof}
$$
\sum_{F\subseteq P}(-1)^{d-\dim(F)}f_k(F)
=\sum_{F\subseteq P}\left(\sum_{\substack{G\subseteq F\\ \dim(G)=k}}(-1)^{d-\dim(F)}\right)
$$
$$
=\sum_{\substack{G\subseteq P\\ \dim(G)=k}}\left(\sum_{F\supseteq G}(-1)^{d-\dim(F)}\right)=
\sum_{\substack{G\subseteq P\\ \dim(G)=k}}\left(\sum_{j=k}^d(-1)^{d-j}f_j([G,P])\right),
$$
where $[G,P]$ denotes the set of faces of $P$ which contain $G$.  By
Euler's equation, the internal sum is equal to zero.
\end{proof}
\end{lem}

We now prove that Theorem~\ref{thm:delta} holds for the $f$-vector of
the summands. The proof follows the same logic, but is slightly
complicated by the fact we are summing on the faces of the sum.

\begin{lem}\label{lem:summand}
Let $P_1,\ldots,P_r$ be $d$-dimensional polytopes and
$P=P_1+\cdots+P_r$ their Minkowski sum. For any face $F$ of $P$, we
denote as $t_i(F)$ be the face of $P_i$ in the decomposition of
$F$. Then for any $1\leq i\leq r$ and $k<d$,
$$
\sum_{F\subseteq P}(-1)^{d-\dim(F)}f_k(t_i(F))=0.
$$
\begin{proof}
$$
\sum_{F\subseteq P}(-1)^{d-\dim(F)}f_k(t_i(F))
=\sum_{F\subseteq P}\left(\sum_{\substack{G\subseteq t_i(F)\\\dim(G)=k}}(-1)^{d-\dim(F)}\right).
$$
For any $F$ and any $i$, $F_i$ is in the decomposition of $F$ if and
only if $\mathcal{N}(F;P)\subseteq \mathcal{N}(F_i;P_i)$. Furthermore,
$G\subseteq F_i$ if and only if $\mathcal{N}(F_i;P_i)\subseteq
cl(\mathcal{N}(G;P_i))$. Therefore, the above expression is equal to
$$
\sum_{F\subseteq P}\left(\sum_{\substack{\cl(\mathcal{N}(G;P_i))\supseteq\mathcal{N}(F;P)\\\dim(G)=k}}(-1)^{d-\dim(F)}\right)
$$
$$
=\sum_{\substack{G\subseteq P_i\\ \dim(G)=k}}\left(\sum_{\mathcal{N}(F;P)\subseteq\cl(\mathcal{N}(G;P_i))}(-1)^{d-\dim(F)}\right).
$$
It is important to note that the sum inside the parentheses is on the
set of faces of $P$ which have their normal cones in the closure of
the normal cone of $G$, which is a face of $P_i$. Here, the polyhedral
cone $\cl(\mathcal{N}(G;P_i))$ is subdivided into a polyhedral complex,
which we call $\Omega(G)$.  We can now write the expression as
$$
\sum_{\substack{G\subseteq P_i\\ \dim(G)=k}}\left(\sum_{j=0}^{d-k}(-1)^j f_j(\Omega(G))\right)=\sum_{\substack{G\subseteq P_i\\ \dim(G)=k}}\chi(\Omega(G)),
$$
where $\chi(\Omega(G))$ is the Euler characteristic of
$\Omega(G)$. The Euler characteristic of any support being independent
of the actual subdivison, we can replace $\Omega(G)$ by its support,
which is $\cl(\mathcal{N}(G;P_i))$:
$$
=\sum_{\substack{G\subseteq P_i\\ \dim(G)=k}}\chi(\cl(\mathcal{N}(G;P_i)))
=\sum_{\substack{G\subseteq P_i\\ \dim(G)=k}}\left(\sum_{j=0}^{d-k}(-1)^j f_j(\cl(\mathcal{N}(G;P_i)))\right)
$$
$$
=\sum_{\substack{G\subseteq P_i\\ \dim(G)=k}}\left(\sum_{j=k}^{d}(-1)^{d-j}f_j([G,P_i]))\right).
$$
Since $\dim(G)=k<d=\dim(P_i)$, the internal sum is equal to zero.
\end{proof}
\end{lem}
Theorem~\ref{thm:delta} is now proved.  We can now prove the main
theorem:
\begin{quote}
\textbf{Theorem~\ref{mainthm}}~
\emph{Let $P_1,\ldots,P_r$ be
$d$-dimensional polytopes relatively in general position, and
$P=P_1+\cdots+P_r$ their Minkowski sum. Then
$$
\sum_{k=0}^{d-1}(-1)^k kf^\delta_k(P) = 0.
$$
}
\end{quote}
\begin{proof}
By Theorem~\ref{thm:delta}, we can replace $f^\delta_k(P)$ by a sum
over proper subfaces:
$$
\sum_{k=0}^{d-1}(-1)^k kf^\delta_k(P) =
\sum_{k=0}^{d-1}(-1)^k k \left(-\sum_{F\subset P}(-1)^{d-\dim(F)}f^\delta_k(F)\right)
$$
$$
= -\sum_{F\subset P}(-1)^{d-\dim(F)}\left(\sum_{k=0}^{d-1}(-1)^k k f^\delta_k(F)\right).
$$
By Lemma~\ref{lemface}, the internal sum is equal to zero.
\end{proof}

We have now proved the main theorem, which makes the assumption that
the summands are full-dimensional. We now extend the result to the
general case. For this, we write an extension of
Lemma~\ref{lem:summand}:
\begin{lem}\label{lem:summandext}
Let $P_1,\ldots,P_r$ be polytopes and $P=P_1+\cdots+P_r$ their
$d$-dimensional Minkowski sum. For any face $F$ of $P$, we denote as
$t_i(F)$ be the face of $P_i$ in the decomposition of $F$. Then for
any $1\leq i\leq r$ and $k<d$,
$$
\sum_{F\subseteq P}(-1)^{d-\dim(F)}f_k(t_i(F))=
(-1)^{d-\dim(P_i)}\delta_{k,\dim(P_i)}.
$$
\begin{proof}
As in the proof of Lemma~\ref{lem:summand}, we prove that
$$ \sum_{F\subseteq P}(-1)^{d-\dim(F)}f_k(t_i(F))
=\sum_{\substack{G\subseteq P_i\\
\dim(G)=k}}\left(\sum_{j=k}^{d}(-1)^{d-j}f_j([G,P_i]))\right).
$$
If $k<\dim(P_i)$, the internal sum is zero as before. If
$k=\dim(P_i)$, then the sums reduce to the single term where $G=P_i$,
and the result is $(-1)^{d-\dim(P_i)}$. If $k>\dim(P_i)$, then the sum
contain no terms.
\end{proof}
\end{lem}
The equation of Theorem~\ref{thm:delta} reads now as
$$
f^\delta_k(P)=-\sum_{F\subset
P}(-1)^{d-\dim(F)}f^\delta_k(F)-(-1)^{d-k}\left|\{P_i\;:\;\dim(P_i)=k\}\right|.
$$
We can now prove the extension:
\begin{quote}
\textbf{Theorem~\ref{maincor2}}~
\emph{Let $P_1,\ldots,P_r$ be polytopes
relatively in general position, and $P=P_1+\cdots+P_r$ their
$d$-dimensional Minkowski sum. Furthermore, let
$S\subseteq\{1,\dots,r\}$ be the set of indices $i$ for which
$\dim(P_i)<d$. Then
$$
\sum_{k=0}^{d-1}(-1)^k k f^\delta_k(P) = (-1)^{d+1}\sum_{i\in S}\dim(P_i).
$$
}
\end{quote}
\begin{proof}
Following the same lines as the proof of Theorem~\ref{mainthm}, we get:
$$
\sum_{k=0}^{d-1}(-1)^k kf^\delta_k(P) =
\sum_{k=0}^{d-1}(-1)^k k
\left(-\sum_{F\subset P}(-1)^{d-\dim(F)}f^\delta_k(F)-(-1)^{d-k}\left|\{P_i\;:\;\dim(P_i)=k\}\right|\right)
$$
$$
= -\sum_{F\subset P}(-1)^{d-\dim(F)}
\underbrace{\left(\sum_{k=0}^{d-1}(-1)^k k f^\delta_k(F)\right)}_0
-\sum_{k=0}^{d-1}(-1)^k k(-1)^{d-k}\left|\{P_i\;:\;\dim(P_i)=k\}\right|
$$
$$
= (-1)^{d+1}\sum_{i\in S}\dim(P_i).
$$
\end{proof}

\section{Application to perfectly centered polytopes}\label{sec:ext}
\begin{thm}\label{fvthm}
Let $P$ be a perfectly centered polytope and $f$ its extended $f$-vector. Then the $f$-vector of
$P+P^*$ can be written as:
$$
f_k(P+P^*) = \sum_{i=0}^kf_{i,i+d-1-k},\quad \forall k=0,\ldots,d-1.
$$
\begin{proof}
Let $P$ be a perfectly centered polytope. From
Theorem~\ref{thm:oldthm}, we know that for every $k$, the number of
$k$-faces of $P+P^*$ is equal to the number of pairs of faces $(G,F)$
of $P$, $G\subseteq F$ so that $\dim(G)+\dim(F^D)=k$, which means
$\dim(G)+d-1-k=\dim(F)$. This is the number of chains of two nontrivial
faces of dimensions $i$ and $i+d-1-k$.
\end{proof}
\end{thm}
We can apply the main theorem to perfectly centered polytopes, which
proves that if $f$ is the $f$-vector of a perfectly centered polytope,
then
$$
\sum_{k=0}^{d-1} (-1)^k k \left(\sum_{i=0}^kf_{i,i+d-1-k}-f_k-f_{d-1-k}\right) = 0.
$$
It turns out that we can extend this result to polytopes and Eulerian
posets in general.

First, let us introduce the Bayer-Billera relations for extended
$f$-vectors:

\begin{lem}[\cite{Bayer85},\cite{Lindstrom71}]\label{dsr}
Let $P$ be an Eulerian poset of rank $d$, $S\subset \{0,\ldots,d-1\}$,
$\{i,k\}\subseteq S \cup \{-1,d\}$, $i<k-1$, and $S$ contains no $j$
so that $i<j<k$. Then

$$
\sum_{j=i+1}^{k-1} (-1)^{j-i-1}f_{S\cup j}(P)=f_S(P)(1-(-1)^{k-i-1}).
$$
\end{lem}
These relations have been presented as an extension of the
Dehn-Sommerville relations. If we examine the special case where
$S=\{i,k\} \subseteq \{-1,\ldots,d\}$, we can write the following
equation:
\begin{lem}\label{dsrshort}
$$
\sum_{j=i}^{k} (-1)^{j}f_{i,j,k}(P)=0.
$$
\end{lem}
Now for the theorem:

\begin{quote}
\textbf{Theorem~\ref{nestthm}}~
\emph{Let $f$ be the extended $f$-vector of an Eulerian poset of rank $d$. Then
$$
\sum_{k=0}^{d-1} (-1)^k k
\left(\sum_{i=0}^kf_{i,i+d-1-k}-f_k-f_{d-1-k}\right) = 0.
$$
}
\end{quote}
\begin{proof}
  Let $P$ be an Eulerian poset of rank $d$. Below, we evaluate and
  rewrite each of the three terms in the parentheses multiplied by
  $(-1)^{(d-1-k)}k$.  (The change of exponent simplifies
  computations):
\begin{equation}\label{eq:first}
-\sum_{k=0}^{d-1}(-1)^{(d-1-k)}k f_k(P) =
\sum_{i=0}^{d-1}(-1)^{d+i}i f_{i,d}(P).
\end{equation}
Then, by using $k'=d-1-k$:
\begin{equation}\label{eq:second}
-\sum_{k=0}^{d-1}(-1)^{(d-1-k)}k f_{d-1-k}(P) =
\sum_{k'=0}^{d-1}(-1)^{k'-1}(d-1-k') f_{-1,k'}(P).
\end{equation}
Finally, by using $k'=i+d-1-k$:
$$
\sum_{k=0}^{d-1}\sum_{i=0}^{k}(-1)^{(d-1-k)}k  f_{i,i+d-1-k}(P)
=\sum_{i=0}^{d-1}\sum_{k'=i}^{d-1}(-1)^{k'-i}(i+d-1-k') f_{i,k'}(P)
$$
$$
=\sum_{i=0}^{d-1}\sum_{k'=i}^{d-1}\left((-1)^{k'+i}i f_{i,k'}(P)+(-1)^{k'+i}(d-1-k') f_{i,k'}(P)\right)
$$
\begin{equation}\label{eq:third}
=\sum_{i=0}^{d-1}\sum_{k'=i}^{d-1}(-1)^{k'+i}i f_{i,k'}(P)+\sum_{k'=0}^{d-1}\sum_{i=0}^{k'}(-1)^{k'+i}(d-1-k') f_{i,k'}(P).
\end{equation}
Combining Equations \ref{eq:first}, \ref{eq:second} and \ref{eq:third}, we get:
$$
\sum_{k=0}^{d-1} (-1)^{(d-1-k)} k \left(\sum_{i=0}^kf_{i,i+d-1-k}-f_k-f_{d-1-k}\right) =
$$
$$
\sum_{i=0}^{d-1}\sum_{k=i}^{d}(-1)^{k+i}i f_{i,k}(P)+\sum_{k'=0}^{d-1}\sum_{i=-1}^{k}(-1)^{k'+i}(d-1-k') f_{i,k}(P)
$$
$$
=\sum_{i=0}^{d-1}i\sum_{k=i}^{d}(-1)^{k+i}f_{i,k,d}(P)+\sum_{k=0}^{d-1}(d-1-k)\sum_{i=-1}^{k}(-1)^{k+i} f_{-1,i,k}(P)=0
$$
The internal sums are zero by the Bayer-Billera relations (Lemma \ref{dsrshort}).
\end{proof}

So we see that this linear relation is a consequence of the
Bayer-Billera relations. This is not surprising, in view of the
theorem of their authors stating that all linear equalities holding
for the extended $f$-vector of Eulerian posets are derived from these equalities
(\cite{Bayer01}).

\begin{ack}
We would like to thank G\"unter M. Ziegler for formulating
Theorem~\ref{maincor2}. We are also grateful to Peter Gritzmann who 
suggested such an extension from Theorem~\ref{mainthm} in the first
place.
\end{ack}

\bibliographystyle{abbrv}

\bibliography{biblio}
\end{document}